\documentclass{article}
\usepackage{amsmath}
\usepackage{amsfonts}
\usepackage{amssymb}
\usepackage{url}
\usepackage{color}
\usepackage{amssymb}
\usepackage{enumerate}

\newtheorem{theorem}{Theorem}
\newtheorem{lemma}[theorem]{Lemma}

\newtheorem{corollary}[theorem]{Corollary}

\newcommand{\Z}{{\mathbb Z}}

\newcommand{\R}{\mathbb R}
\newcommand{\FourTiTwo}{{\tt 4ti2}}

\newcommand{\Graver}{{{\cal G}}}
\newcommand{\UGB}{{{\cal U}}}

\newcommand{\Groebner}{{{\cal G}_\prec}}
\newcommand{\Orthant}{\mathbb O}
\newcommand{\Ideal}{{{\cal I}}}
\newcommand{\boproof}{{\bf Proof.} }
\newcommand{\eoproof}{\hspace*{\fill} $\square$ \vspace{5pt}}

\DeclareMathOperator{\conv}{conv}
\DeclareMathOperator{\supp}{supp}

\usepackage{enumerate}

\begin{document}
\setlength{\parindent}{0pt} \setlength{\parskip}{2ex plus 0.4ex
minus 0.4ex}

\title{On the Gr\"obner complexity of matrices}

\author{{\bf Raymond~Hemmecke}\\University of Magdeburg, Germany\\
\\
{\bf Kristen~A.~Nairn}\\College of St. Benedict, MN, USA}

\date{}

\maketitle

\begin{abstract}
In this paper we show that if for an integer matrix $A$
the universal Gr\"obner basis of the associated toric ideal $\Ideal_A$
coincides with the Graver basis of $A$, then the Gr\"obner complexity
$u(A)$ and  the Graver complexity $g(A)$ of its higher Lawrence
liftings agree, too. We conclude that for the matrices $A_{3\times 3}$
and  $A_{3\times 4}$, defining the $3\times 3$ and $3\times 4$
transportation problems, we have $u(A_{3\times 3})=g(A_{3\times 3})=9$
and $u(A_{3\times 4})=g(A_{3\times 4})\geq 27$. Moreover, we prove
$u(A_{a,b})=g(A_{a,b})=2(a+b)/\gcd(a,b)$ for positive 
integers $a,b$ and $A_{a,b}=\left(\begin{smallmatrix}1 & 1 & 1 &
1\\0 & a & b & a+b\end{smallmatrix}\right)$.  

\end{abstract}

\section{Introduction}
\label{Section: Introduction}

In this paper we deal with Graver bases and universal Gr\"obner bases
associated to a matrix. The Graver basis of $A\in\Z^{d\times n}$ is
defined as the union 
\[
\Graver(A):=\bigcup_{j=1}^{2^n} H_j\setminus\{0\}
\]
of the inclusion-minimal Hilbert bases $H_j$ of the pointed rational
polyhedral cones 
\[
C_j:=\ker(A)\cap\Orthant_j=\{z\in\Orthant_j:Az=0\}
\]
as $\Orthant_j$ ranges over all $2^n$ orthants of $\R^n$
\cite{Graver:75}. Moreover, we call
\[
\Ideal_A:=\langle x^u-x^v:Au=Av,u,v\in\Z^n_+\rangle
\]
the \emph{toric ideal} associated to $A$, and for a given term
ordering $\prec$, we call $\Groebner(A)$ a minimal reduced Gr\"obner
basis of $A$ with respect to $\prec$, if
$\{x^{u^+}-x^{u^-}:u\in\Groebner(A)\}$ is a minimal reduced
Gr\"obner basis of $\Ideal_A$ with respect to $\prec$. By $\UGB(A)$
we denote the universal Gr\"obner basis of $A$, being the union over
all minimal reduced Gr\"obner bases of $A$. Note that one can show
that the relation $\Groebner(A)\subseteq\UGB(A)\subseteq\Graver(A)$
holds for any term ordering $\prec$. In particular, $\UGB(A)$ is 
finite.

In \cite{Santos+Sturmfels}, Santos and Sturmfels dealt with the
question of how complex or complicated the Graver bases of the
matrices 
\[
A^{(N)}:= \left(
\begin{array}{cccc}
  I_n & I_n & \cdots & I_n \\
  A &   &        &   \\
    & A &        &   \\
    &   & \ddots &   \\
    &   &        & A \\
\end{array}
\right)
\]
become as $N\in\Z_+$ grows. For a vector $x=(x^1,\ldots,x^N)$ with
$x^i\in\Z^n$ for $i=1,\ldots,N$, we call $x^1,\ldots,x^N$ the
\emph{layers} of $x$. Moreover, we call the number $|\{i:x^i\neq 0\}|$
of nonzero layers of $x$ the \emph{type} of $x$. With these notions,
Santos and Sturmfels showed that there is a constant $g(A)$, depending
only on $A$ but not on $N$, such that the types of the Graver basis
elements of $A^{(N)}$ are bounded by $g(A)$ for all $N$. They coined 
the notion ``Graver complexity'' for this constant $g(A)$. Moreover, they
presented an algorithm to compute $g(A)$, with which
they computed the Graver complexity of the matrix
$\left(\begin{smallmatrix}1 & 1 & 1 & 1\\0 & 1 & 2&
3\end{smallmatrix}\right)$ defining the twisted cubic, which is $6$, 
and the Graver  complexity of the matrix $A_{3\times 3}$ defining the 
$3\times 3$ transportation polytope, which is $9$. Already the next 
bigger case of $3\times 4$ transportation polytopes was left open. 

In this paper, we consider the analogous notion of \emph{Gr\"obner
complexity} $u(A)$ of $A$ as the maximal type of an element in
$\UGB(A^{(N)})$ for all $N$. By the results of \cite{Aoki+Takemura} and
\cite{Santos+Sturmfels}, we have $5\leq u(A_{3\times 3})\leq 9$. Boffi
and Rossi \cite{Boffi+Rossi} proved that the maximal type of a vector
appearing in any \emph{lexicographic} Gr\"obner basis of $A_{3\times
3}^{(N)}$, $N\geq 5$, is $5$. This left open the question whether
there exist other term orderings such that the corresponding
Gr\"obner bases contain a vector of type $6$, $7$, $8$, or $9$.

In Section \ref{Section: Proof of main theorem}, we prove our main
result of this paper. In fact, with Lemma \ref{Lemma: Main lemma}, we
even prove a deeper structural result on $\UGB(A^{(N)})$, from which
Theorem \ref{Theorem: Main theorem} follows by the results in
\cite{Santos+Sturmfels}. 

\begin{theorem}\label{Theorem: Main theorem}
Let $A\in\Z^{d\times n}$. If the universal Gr\"obner basis $\UGB(A)$
and the Graver basis $\Graver(A)$ coincide, then $u(A)=g(A)$, that is,
Gr\"obner complexity and Graver complexity of $A$ are equal.
\end{theorem}

Note that we do not claim that the universal Gr\"obner bases and the
Graver bases of $A^{(N)}$ are the same for each $N$. In fact, we leave
this as an open question that remains to be clarified. Our theorem has a
few nice consequences. For example, as $\UGB(A)=\Graver(A)$
whenever $A$ is a unimodular matrix, we get the following nice fact. 

\begin{corollary}\label{Corollary: Unimodular matrices}
For unimodular matrices, Gr\"obner complexity and Graver complexity
are equal.  
\end{corollary}

In particular, this implies $u(A_{3\times 3})=g(A_{3\times 3})=9$. In
fact, in Section \ref{Section: 3x3} below, we explicitly state
elements in $\UGB(A_{3\times 3}^{(9)})$ of types $6$, $7$, $8$, and
$9$, together with term orderings for which these elements 
appear in the corresponding Gr\"obner bases. It comes as a little
surprise that there are indeed elements in $\UGB(A_{3\times 3}^{(9)})$
that are more complicated (= have a bigger type) than the  
elements in any lexicographic Gr\"obner basis of $A_{3\times
3}^{(9)}$. 

In Section \ref{Section: 3x4}, we consider the case of
$A_{3\times 4}$ and show the following. Note that $u(A_{3\times
4})=g(A_{3\times 4})$ already follows from Corollary \ref{Corollary:
Unimodular matrices}, as $A_{3\times 4}$ is unimodular.

\begin{corollary}\label{Corollary: 3x4}
For $A_{3\times 4}$, we have $u(A_{3\times 4})=g(A_{3\times 4})\geq 27$.
\end{corollary}

In fact, we conjecture this bound to be tight, that is, $u(A_{3\times
4})=g(A_{3\times 4})=27$. 

Finally, in Sections \ref{Section: A_a,b} and \ref{Section:
2c-conjecture}, we show the following result.  

\begin{lemma}\label{Lemma: A_a,b lemma}
Let $a,b\in\Z_{>0}$ and $A_{a,b}=\left(\begin{smallmatrix}1 & 1 & 1 &
1\\0 & a & b &  a+b\end{smallmatrix}\right)$. Then
$u(A_{a,b})=g(A_{a,b})= 2(a+b)/\gcd(a,b)$. 
\end{lemma}

Consequently, for $a=1$ and $b=2$, we conclude that Gr\"obner
complexity and Graver complexity of the matrix defining the twisted 
cubic both equal $6$. To prove Lemma \ref{Lemma: A_a,b lemma}, 
we first show the inequality $g(A_{a,b})\geq u(A_{a,b})\geq 2(a+b)/\gcd(a,b)$ 
in Section \ref{Section: A_a,b}. Then, in Section 
\ref{Section: 2c-conjecture}, we show that in
fact $g(A_{a,b})=2(a+b)/\gcd(a,b)$, finally settling an open problem
from \cite{Nairn:PhD-thesis}.

\section{Proof of main theorem}\label{Section: Proof of main theorem}  

Let us now prove Theorem \ref{Theorem: Main theorem}. In fact, we show
a stronger result from which the statement of Theorem \ref{Theorem:
Main theorem} follows immediately by the results of
\cite{Santos+Sturmfels}. 

\begin{lemma}\label{Lemma: Main lemma}
Let $A\in\Z^{d\times n}$ and let $\UGB(A)=\{g_1,\ldots,g_k\}$ denote
the elements in the universal Gr\"obner basis of $A$. We assume that
$\UGB(A)$ is symmetric, that is, if $g\in\UGB(A)$ then also
$-g\in\UGB(A)$. 

Let $\lambda\in\Z^k_+$ be the coefficient vector of a minimal 
nonnegative integer relation among $\{g_1,\ldots,g_k\}$, that is,
$\sum_{i=1}^k \lambda_ig_i=0$. With $s=|\supp(\lambda)|$, the vector
$x\in\ker(A^{(s)})$ shall denote an arrangement of $\lambda_1$ layers
$g_1$, $\lambda_2$ layers $g_2$, and so on, in any arbitrary but fixed
order. Then each such vector $x$ belongs to $\UGB(A^{(s)})$. 
\end{lemma}

\boproof
For each $g_i\in\UGB(A)$, let ${\cal F}_i$ denote the polyhedron
$\conv(\{z\in\Z^n_+:Az=A(x^i)^+\})$. For each $g_i\in\UGB(A)$ there
exists some vector $c_i\in\R^n$ and some number $\gamma_i\in\R$ such
that the inequality $c_i^\intercal z\geq\gamma_i$ defines an egde of
${\cal F}_i$ with edge direction $g_i$. In fact ${\cal
F}_i=\conv(\{g_i^+,g_i^-\})$, see \cite{Sturmfels+Thomas:97}.

Now consider the vector $c\in\R^{sn}$ formed out of $\lambda_1$ copies
of $c_1$, $\lambda_2$ copies of $c_2$, and so on, in the same order of
indices as $x$ was formed. Then, by construction, $c^\intercal
z\geq\sum_{i=1}^k\lambda_i\gamma_i$ is a valid inequality of the
polyhedron ${\cal P}=\conv(\{z\in\Z^{sn}_+:A^{(s)}z=A^{(s)}x^+\})$ and
hence defines a face ${\cal F}$ of it. Again by construction, any
lattice point on this face ${\cal F}$ can only have $\lambda_i$ layers
(in total) of $g_i^+$ and  $g_i^-$ for each $i=1,\ldots, k$, as otherwise
any strict (face defining) inequality for some layer would imply the
relation $c^\intercal z>\sum_{i=1}^k\lambda_i\gamma_i$ for the whole
vector. 
 
Now assume that for a lattice point $y$ of ${\cal F}$, we choose
$\mu_1$ layers $g_1^+$ and $\lambda_1-\mu_1$ layers $g_1^-$,
$\mu_2$ layers $g_2^+$ and $\lambda_2-\mu_2$ layers $g_2^-$, and
so on, in the same order of indices as in $x$. Thus, as
$y\in{\cal F}$, we have $\sum_{i=1}^k 
\mu_ig_i^++(\lambda_i-\mu_i)g_i^-=A^{(s)}x^+=A^{(s)}x^-=\sum_{i=1}^k
\lambda_ig_i^-$. Consequently, we get $0=\sum_{i=1}^k \mu_i
(g_i^+-g_i^-)=\sum_{i=1}^k \mu_ig_i$. This is a contradiction 
to the minimality of $\lambda$ unless $\mu=0$ or
$\mu=\lambda$. Consequently, ${\cal F}$ contains only two lattice
points, namely $x^+$ (for $\mu=\lambda$) and $x^-$ (for
$\mu=0$). Thus, ${\cal F}$ is an edge of ${\cal P}$ with edge
direction $x$. Therefore, $x$ belongs to $\UGB(A^{(s)})$, see
\cite{Sturmfels+Thomas:97}. \eoproof

Let now $\UGB(A)=\Graver(A)$. Then, by the results of
\cite{Santos+Sturmfels}, any element in $\Graver(A^{(g(A))})$ of
maximal type corresponds to a vector $x$ as described in Lemma
\ref{Lemma: Main lemma}. Thus, Theorem \ref{Theorem: Main theorem} is
proved, too.

\section{Elements in $\UGB\left(A_{3\times3}^{(9)}\right)$ of types
$6,7,8$, and $9$}\label{Section: 3x3}

In this section, we present elements
$x_6,x_7,x_8,x_9\in\UGB\left(A_{3\times3}^{(9)}\right)$ of types
$6,7,8$, and $9$. The nonzero layers of the four elements
$x_6,x_7,x_8,x_9\in\UGB\left(A_{3\times3}^{(9)}\right)$ are

{\tiny
\begin{eqnarray*}
z_6 & : & \left(
\begin{array}{rrr}
0 &  1 & -1\\
0 &  0 & 0\\
0 & -1 & 1\\
\end{array}
\right) +\left(
\begin{array}{rrr}
0 & 1 & -1\\
1 & -1 & 0\\
-1 & 0 & 1\\
\end{array}
\right) +\left(
\begin{array}{rrr}
1 & -1 & 0\\
-1 & 0 & 1\\
0 & 1 & -1\\
\end{array}
\right)+\left(
\begin{array}{rrr}
1 & -1 & 0\\
0 & 1 & -1\\
-1 & 0 & 1\\
\end{array}
\right)+2\cdot\left(
\begin{array}{rrr}
-1 & 0 & 1\\
0 & 0 & 0\\
1 & 0 & -1\\
\end{array}
\right),\\
z_7 & : & \left(
\begin{array}{rrr}
0 & 1 & -1\\
1 & -1 & 0\\
-1 & 0 & 1\\
\end{array}
\right) +\left(
\begin{array}{rrr}
1 & -1 & 0\\
0 & 1 & -1\\
-1 & 0 & 1\\
\end{array}
\right) +\left(
\begin{array}{rrr}
1 & 0 & -1\\
-1 & 1 & 0\\
0 & -1 & 1\\
\end{array}
\right)+\left(
\begin{array}{rrr}
1 & 0 & -1\\
0 & -1 & 1\\
-1 & 1 & 0\\
\end{array}
\right)+3\cdot\left(
\begin{array}{rrr}
1 & 0 & -1\\
0 & -1 & 1\\
-1 & 1 & 0\\
\end{array}
\right),\\
z_8 & : & \left(
\begin{array}{rrr}
1 & -1 & 0\\
-1 & 0 & 1\\
0 & 1 & -1\\
\end{array}
\right) +\left(
\begin{array}{rrr}
1 & 0 & -1\\
0 & -1 & 1\\
-1 & 1 & 0\\
\end{array}
\right) +\left(
\begin{array}{rrr}
-1 & 0 & 1\\
1 & -1 & 0\\
0 & 1 & -1\\
\end{array}
\right)+2\cdot\left(
\begin{array}{rrr}
1 & -1 & 0\\
0 & 1 & -1\\
-1 & 0 & 1\\
\end{array}
\right)+3\cdot\left(
\begin{array}{rrr}
-1 & 1 & 0\\
0 & 0 & 0\\
1 & -1 & 0\\
\end{array}
\right),\\
z_9 & : & \left(
\begin{array}{rrr}
-1 & 0 & 1\\
0 & 1 & -1\\
1 & -1 & 0\\
\end{array}
\right) +\left(
\begin{array}{rrr}
0 & 1 & -1\\
-1 & 0 & 1\\
1 & -1 & 0\\
\end{array}
\right) +2\cdot\left(
\begin{array}{rrr}
-1 & 0 & 1\\
1 & -1 & 0\\
0 & 1 & -1\\
\end{array}
\right)+2\cdot\left(
\begin{array}{rrr}
0 & 1 & -1\\
1 & -1 & 0\\
-1 & 0 & 1\\
\end{array}
\right) +3\cdot\left(
\begin{array}{rrr}
1 & -1 & 0\\
-1 & 1 & 0\\
0 & 0 & 0\\
\end{array}
\right).
\end{eqnarray*}
}

Note that due to the underlying symmetry of the problem matrix
$A_{3\times 3}^{(9)}$, any arrangement of these $6$, $7$, $8$, or
$9$ layers together with sufficiently many zero layers 
gives an element in $\UGB\left(A_{3\times3}^{(9)}\right)$. Fix
any such arrangement, for example first using the first layer type,
then the second layer type, and so on, and call the resulting vectors
$x_6$, $x_7$, $x_8$, and $x_9$. 
By Lemma \ref{Lemma: Main lemma}, $x_6$, $x_7$, $x_8$, and $x_9$
belong to $\UGB\left(A_{3\times 3}^{(9)}\right)$. Valid
inequalities defining the edges
$\conv\left(\left\{x_i^+,x_i^-\right\}\right)$ of ${\cal F}_j:= 
\conv\left\{z\in\Z^{81}_+:A_{3\times3}^{(9)}\;z=A_{3\times3}^{(9)}\;
x_j^+\right\}$, $j=6,7,8,9$, are given by
$\sum_{i\in\{1,\ldots,81\}\setminus\supp(x_j)} z_i\geq 0$.

One may now ask what term ordering one has to choose to obtain
$x_6$, $x_7$, $x_8$, or $x_9$ as a Gr\"obner basis element. For
this, one may use
$\sum_{i\in\{1,\ldots,81\}\setminus\supp(x_j)} e_i$, $j=6,7,8,9$, as
cost vector and any term ordering to break ties. Herein, $e_i$ denotes 
as usual the $i$-th unit vector. For example, for
$j=9$, using the ${\tt groebner}$ function of \FourTiTwo{} (version
1.3.1) and the default degrevlex ordering (of \FourTiTwo) for
tie-breaking, one obtains the $218,785$ vectors in the corresponding
Gr\"obner basis within $51$ minutes on an AMD Opteron 2.4 GHz CPU
running SuSE Linux 10.0.

\section{Proof of $u(A_{3\times 4})=g(A_{3\times 4})\geq
27$}\label{Section: 3x4}   

In this section, we present an element
$x_{27}\in\UGB\left(A_{3\times4}^{(27)}\right)$ and thus prove that
$u(A_{3\times 4})=g(A_{3\times 4})\geq 27$. The nonzero layers of
$x_{27}\in\UGB\left(A_{3\times 4}^{(27)}\right)$ are 

{\tiny
\begin{eqnarray*}
z_{27} & : & 1\cdot\left(
\begin{array}{rrr}
0 & -1 & 1\\
0 & 0 & 0\\
1 & 0 & -1\\
-1 & 1 & 0\\
\end{array}
\right) +2\cdot\left(
\begin{array}{rrr}
-1 & 1 & 0\\
1 & 0 & -1\\
0 & -1 & 1\\
0 & 0 & 0\\
\end{array}
\right) +3\cdot\left(
\begin{array}{rrr}
1 & -1 & 0\\
0 & 1 & -1\\
0 & 0 & 0\\
-1 & 0 & 1\\
\end{array}
\right)+3\cdot\left(
\begin{array}{rrr}
0 & 1 & -1\\
1 & -1 & 0\\
0 & 0 & 0\\
-1 & 0 & 1\\
\end{array}
\right)+\\&&5\cdot\left(
\begin{array}{rrr}
0 & 1 & -1\\
-1 & 0 & 1\\
1 & -1 & 0\\
0 & 0 & 0\\
\end{array}
\right)+6\cdot\left(
\begin{array}{rrr}
1 & -1 & 0\\
0 & 0 & 0\\
-1 & 0 & 1\\
0 & 1 & -1\\
\end{array}
\right)+7\cdot\left(
\begin{array}{rrr}
-1 & 0 & 1\\
0 & 0 & 0\\
0 & 1 & -1\\
1 & -1 & 0\\
\end{array}
\right).
\end{eqnarray*}
}

Again, due to symmetry, the actual arrangement of the $27$ layers is 
not important.
Thus, we may again assume that we first use the first layer type,
then the second layer type, and so on, and call the resulting vector
$x_{27}$. To show that indeed $x_{27}\in\UGB\left(A_{3\times
4}^{(27)}\right)$, one only has to check that $x_{27}$ is given by a
minimal relation among the elements in $\UGB(A_{3\times
4})=\Graver(A_{3\times 4})$ as required by Lemma \ref{Lemma: Main lemma}. 
This is a feasibility problem in only $7$ integer variables, which can
easily be solved using the {\tt zsolve} function of \FourTiTwo{} or,
with a bit more work, even by hand.

\section{Proof of $g(A_{a,b})\geq u(A_{a,b})\geq 2(a+b)$}   
\label{Section: A_a,b}    

In this section, we prove $g(A_{a,b})\geq u(A_{a,b})\geq
2(a+b)/\gcd(a,b)$. In the next section, we show that in fact equality
holds, implying Lemma \ref{Lemma: A_a,b lemma}. Let us remind the
reader that $1\leq a<b$ are positive integers and we 
consider the matrix 
$A_{a,b}=\left(\begin{smallmatrix}1 & 1 & 1 & 1\\0 & a & b &
a+b\end{smallmatrix}\right)$. Note that we may divide the second row
of $A_{a,b}$ by $\gcd(a,b)$ without changing the integer kernel of the
matrix. Thus, we may without loss of generality assume $a$ and $b$ to be
coprime. In order to prove $g(A_{a,b})\geq u(A_{a,b})\geq 2(a+b)$
using Lemma \ref{Lemma: Main lemma}, we first prove that
$(b-1,-a-b+1,1,a-1)^\intercal$, $(-b,a+b,0,-a)^\intercal$, and
$(a,0,-a-b,b)^\intercal$ belong to $\UGB(A_{a,b})$.  

Let us first consider $x=(-b,a+b,0,-a)^\intercal$ and the face of
$\conv(\{y\in\Z^4_+:A_{a,b}y=A_{a,b}x^+\})$ defined by the valid
inequality $y_3\geq 0$. The linear system defining this face reads
\begin{eqnarray*}
y_1+y_2+y_3+y_4 & = & a+b\\
ay_2+by_3+(a+b)y_4 & = & a(a+b)\\
y_3 & = & 0,
\end{eqnarray*}
with nonnegative integers $y_1,\ldots,y_4$. Eliminating $y_3$ and
subtracting $a$ times the first equation from the second equation,
we obtain the equivalent system
\begin{eqnarray*}
y_1+y_2+y_4 & = & a+b\\
-ay_1+by_4 & = & 0\\
y_3 & = & 0.
\end{eqnarray*}
As $\gcd(a,b)=1$ and since $0\leq y_1,y_4\leq a+b$, by the first
equation, we conclude that the second equation has only two
solutions $y_1=0, y_4=0$ and $y_1=b, y_4=a$. In the first case, we
obtain $y_2=a+b$ and in the second case $y_2=0$. Thus, the face under
consideration is $\conv(\{x^+,x^-\})$ and hence $x\in\UGB(A_{a,b})$.

For $x=(a,0,-a-b,b)^\intercal$, the arguments are similar (due to the 
symmetry in $a$ and $b$). Hence again $x\in\UGB(A_{a,b})$.

For $x=(b-1,-a-b+1,1,a-1)^\intercal$, the proof is a bit more
complicated. Let us consider the face of
$\conv(\{y\in\Z^4_+:A_{a,b}y=A_{a,b}x^+\})$ defined by the valid 
inequality $(a-1)y_3-y_4\geq 0$. Below, we will see that
$(a-1)y_3-y_4\geq 0$ is indeed valid. The defining linear systems
reads 
\begin{eqnarray*}
y_1+y_2+y_3+y_4 & = & a+b-1\\
ay_2+by_3+(a+b)y_4 & = & a(a+b-1),
\end{eqnarray*}
with nonnegative integers $y_1,\ldots,y_4$. Subtracting $a$ times
the first equation from the second equation, we obtain the
equivalent system
\begin{eqnarray*}
y_1+y_2+y_3+y_4 & = & a+b-1\\
-ay_1+(b-a)y_3+by_4 & = & 0.
\end{eqnarray*}
As $\gcd(a,b)=1$, we conclude from the second equation that
$a(y_1+y_3)$ and thus also $y_1+y_3$ is divisible by $b$. Since
$a<b$, we conclude from the first equation that either $y_1+y_3=0$
or $y_1+y_3=b$. In the first case, we obtain $y_1=y_3=y_4=0, y_2=a+b-1$, 
and in the second case, we easily get
$(b-a,a-1,a,0)^\intercal,(b-a+1,a-2,a-1,1)^\intercal,\ldots,
(b-1,0,1,a-1)^\intercal$ as the solutions for
$(y_1,\ldots,y_4)^\intercal$. This proves that the inequality
$(a-1)z_3-z_4\geq 0$ is indeed valid for the fiber
$\conv(\{y\in\Z^4_+:A_{a,b}y=A_{a,b}x^+\})$. Moreover, our discussion
shows that it defines the face $\conv(\{x^+,x^-\})$. Therefore, we
conclude $x\in\UGB(A_{a,b})$.

Finally, to prove $g(A_{a,b})\geq u(A_{a,b})\geq 2(a+b)$ using Lemma
\ref{Lemma: Main lemma}, let us consider the vector $x$ from
$\ker\left(A_{a,b}^{(2a+2b)}\right)$ given by the following $2(a+b)$
layers
\[
(a+b)\cdot\left(
\begin{array}{c}
b-1\\
-a-b+1\\
1\\
a-1\\
\end{array}
\right)+(a+b-1)\cdot\left(
\begin{array}{c}
-b\\
a+b\\
0\\
-a\\
\end{array}
\right)+1\cdot\left(
\begin{array}{c}
a\\
0\\
-a-b\\
b\\
\end{array}
\right).
\]
Again, the actual arrangement of these layers is not important due
to the symmetry underlying $A_{a,b}^{(2a+2b)}$. 

One easily checks that this relation of three elements in
$\UGB(A_{a,b})$ is indeed minimal. Clearly, there is no relation among
only two of the elements. Thus, $(a,0,-a-b,b)^\intercal$ has to be
used once and therefore, looking at the third coordinate,
$(b-1,-a-b+1,1,a-1)^\intercal$ has to be used $a+b$
times. Consequently, the coefficient of the vector
$(-b,a+b,0,-a)^\intercal$ is $a+b-1$, and we obtain the initial
relation which is thus minimal. Therefore, by Lemma
\ref{Lemma: Main lemma}, $x\in\UGB\left(A_{a,b}^{(2a+2b)}\right)$ and
hence $g(A_{a,b})\geq u(A_{a,b})\geq 2(a+b)$.

\section{Proof of $g(A_{a,b})=2(a+b)$}\label{Section: 2c-conjecture}   

In this section we show the following.

\begin{lemma}\label{Lemma: g=2a+2b}
If $1\leq a<b$ are coprime integers, then $g(A_{a,b})=2(a+b)$. 
\end{lemma}

In order to compute $g(A_{a,b})$ via the construction from
\cite{Santos+Sturmfels}, we need to first compute the Graver basis of
$A_{a,b}$ and write down its elements as the columns of a new matrix
$G_{a,b}$. (Note that, as $\Graver(A_{a,b})$ is symmetric, we only
need to choose one vector out of each pair
$g,-g\in\Graver(A_{a,b})$. See \cite{Santos+Sturmfels} for more
details.) Then $g(A_{a,b})$ is equal to the maximum $1$-norm of the
vectors appearing in the Graver basis $\Graver(G_{a,b})$. Let us start
by presenting the Graver basis of $A_{a,b}$. 

\begin{lemma}
Let $v=(-b,a+b,0,-a)^\intercal$ and $h=(1,-1,-1,1)^\intercal$. Then we
have 
\[
\Graver(A_{a,b})=\pm\{v,v+h,v+2h,\ldots,v+(a+b)h,h\}.
\]
\end{lemma}

\boproof To show our claim, we only need to check the
criteria of Lemma 2 in \cite{Hemmecke:PSP}: $v$ and $h$ do indeed
generate $\ker_Z(A_{a,b})$ over $\Z$, the set $\Graver(A_{a,b})$ is
indeed symmetric (that is $g\in\Graver(A_{a,b})$ implies
$-g\in\Graver(A_{a,b})$), and for each choice of
$g_1,g_2\in\Graver(A_{a,b})$, the vector $g_1+g_2$ can be 
written as a sign-compatible positive integer linear combination of
elements in $\Graver(A_{a,b})$. Due to the simple structure of
$\Graver(A_{a,b})$, the latter requires only an easy case
study on the possible sign patterns of $g_1+g_2$. \eoproof

Now we need to find the maximum $1$-norm among the vectors in the
Graver basis of the matrix $G_{a,b}=(v,v+h,v+2h,\ldots,v+(a+b)h,h)$. The
following lemma tells us that we may consider the matrix
$B_{a+b}=\left(\begin{smallmatrix}1 & 1 & 1 & \ldots & 1 & 0\\0 & 1 &
2 & \ldots & a+b & 1\end{smallmatrix}\right)$ instead.  

\begin{lemma}
If $1\leq a<b$ are coprime integers, then
$\ker_\Z(G_{a,b})=\ker_\Z(B_{a+b})$. 
\end{lemma}

\boproof
First note that $G_{a,b}=(v,h)^\intercal B_{a+b}$. Now let
$x\in\Z^{a+b+2}$. Since $v$ and $h$ are linearly independent
(they generate $\ker_Z(A_{a,b})$ over $\Z$), we conclude
\[
x\in\ker_\Z(G_{a,b})\Leftrightarrow G_{a,b}x=0
\Leftrightarrow (v,h)^\intercal [B_{a+b}x]=0 \Leftrightarrow
B_{a+b}x=0\Leftrightarrow x\in\ker_\Z(B_{a+b}).
\]
\eoproof

It remains to show the following fact. Our claim
$g(A_{a,b})=2(a+b)$ will then follow immediately. Note that this lemma
proves the $2c$-conjecture from \cite{Nairn:PhD-thesis}, Corollary
\ref{Corollary: 2c-conjecture}. 

\begin{lemma}\label{Lemma: 2n-2 bound}
The maximum $1$-norm of a vector appearing in the Graver basis of the 
matrix $A_n=\left(\begin{smallmatrix}1 & 1 & 1 & \ldots & 1 & 0\\
1 & 2 & 3 & \ldots & n & 1\end{smallmatrix}\right)$ is $2(n-1)$.
\end{lemma}

\boproof
Let $x\in\ker(A_n)$ with $x_{n+1}\geq 0$. Then we can translate $x$
uniquely into a relation among the numbers $1,2,\ldots,n$ of the form
$\sum_{i=1}^ka_i+\sum_{j=1}^l 1 = \sum_{i=1}^kb_i$, with $a_i\neq
b_j$ for $i,j=1,\ldots,k$, where for $t=1,\ldots,n$, $x_t$ counts the
number of occurrences of $a_i=t$ minus the number of occurrences of
$b_i=t$, and where we have $l=x_{n+1}$. Vice-versa, we can translate any
relation $\sum_{i=1}^ka_i+\sum_{j=1}^l 1 = \sum_{i=1}^kb_i$ uniquely back 
into a vector $x\in\ker(A_n)$ by simple counting.

Assume now that $x\in\Graver(A_n)$. This implies that 
there does not exisst a non-trivial sub-identity 
$\sum_{i\in I}a_i+\sum_{j\in J} 1 = \sum_{i\in I}b_i$ in 
$\sum_{i=1}^ka_i+\sum_{j=1}^l 1 = \sum_{i=1}^kb_i$, as otherwise the 
corresponding vector of counts would contradict the minimality of $x$.

Without loss of generality, we may assume that $a_1\leq\ldots\leq a_k$
and $b_1\leq\ldots\leq b_k$. Now define $\delta_i=a_i-b_i$,
$i=1,\ldots,k$. Note that there cannot exist a sub-identity $\sum_{i\in
I}\delta_i+\sum_{j\in J} 1=0$ of $\sum_{i=1}^k\delta_i+\sum_{j=1}^l 1
=0$, since then $\sum_{i\in I}a_i+\sum_{j\in J} 1 = \sum_{i\in
I}b_i$ would contradict the minimality of
$\sum_{i=1}^ka_i+\sum_{j=1}^l 1 = \sum_{i=1}^kb_i$. Thus, by bringing
all negative values $\delta_i$ to the right-hand side of the relation, 
we obtain a \emph{primitive partition identity} $\sum_{i=1,
\delta_i>0}^k\delta_i+\sum_{j=1}^l 1 =\sum_{i=1,
\delta_i<0}^k(-\delta_i)$. By Corollary 1 in
\cite{Diaconis+Graham+Sturmfels}, we obtain for our primitive
partition identity the bound $k+l\leq\Delta_++\Delta_-$, where
$\Delta_+=\max\{\delta_i:\delta_i>0\}\geq 1$ and
$\Delta_-=\max\{-\delta_i:\delta_i<0\}\geq 1$.  

The remainder of our proof follows nearly literally the proof of Theorem
3 in \cite{Diaconis+Graham+Sturmfels}. For the benefit of the reader,
we include these few lines here. Let $i_0$ and $j_0$ be such that
$b_{i_0}-a_{i_0}=\Delta_-$ and $a_{j_0}-b_{j_0}=\Delta_+$. Now we
distinguish two cases. 

If $i_0<j_0$, then
\[
1+\Delta_-\leq a_{i_0}+\Delta_-=b_{i_0}\leq
b_{j_0}=a_{j_0}-\Delta_+\leq n-\Delta_+.
\]
If $i_0>j_0$, then
\[
n-\Delta_-\geq b_{i_0}-\Delta_-=a_{i_0}\geq
a_{j_0}=-b_{j_0}+\Delta_+\geq 1+\Delta_+.
\]
In both cases, we obtain $\Delta_++\Delta_-\leq n-1$. Consequently,
\[
\|x\|_1=2k+l\leq 2(k+l)\leq 2(\Delta_++\Delta_-)\leq 2(n-1).
\]

Finally, we should mention that this upper bound is tight, since
$1\cdot 1-(n-1)\cdot (n-1) + (n-2)\cdot n=0$ corresponds to the 
minimal vector $e_1-(n-1)e_(n-1)+(n-2)e_n$ in $\Graver(A_n)$.  
\eoproof

\begin{corollary}\label{Corollary: 2c-conjecture}
The maximum $1$-norm of a vector appearing in the Graver basis of the 
matrix $\left(\begin{smallmatrix}1 & 1 & 1 & \ldots & 1 & 0\\
0 & 1 & 2 & \ldots & c & 1\end{smallmatrix}\right)$ is $2c$.
\end{corollary}

\boproof
Note that
$
\ker\left(\begin{smallmatrix}1 & 1 & 1 & \ldots & 1 & 0\\
0 & 1 & 2 & \ldots & c & 1\end{smallmatrix}\right)=
\ker\left(\begin{smallmatrix}1 & 1 & 1 & \ldots & 1 & 0\\
1 & 2 & 3 & \ldots & c+1 & 1\end{smallmatrix}\right)
$ and thus
$
\Graver\left(\begin{smallmatrix}1 & 1 & 1 & \ldots & 1 & 0\\
0 & 1 & 2 & \ldots & c & 1\end{smallmatrix}\right)=
\Graver\left(\begin{smallmatrix}1 & 1 & 1 & \ldots & 1 & 0\\
1 & 2 & 3 & \ldots & c+1 & 1\end{smallmatrix}\right).
$
Now apply Lemma \ref{Lemma: 2n-2 bound}. \eoproof

Applying this corollary to our matrix $\Graver(A_{a,b})$, we conclude
that $g(A_{a,b})=2(a+b)$, and Lemma \ref{Lemma: g=2a+2b} is
proved. Consequently, $u(A_{a,b})=g(A_{a,b})=2(a+b)$, as claimed in
Lemma \ref{Lemma: A_a,b lemma}.


\begin{thebibliography}{99}

\bibitem{4ti2}
4ti2 team. {\tt 4ti2}--A software package for algebraic, geometric and
combinatorial problems on linear spaces. Available at www.4ti2.de.

\bibitem{Aoki+Takemura}
S.~Aoki and A.~Takemura.
Minimal basis for connected Markov chain over $3\times3\times K$
contingency tables with fixed two-dimensional marginals.
{\sl Austr. New Zeal. J. Stat.} \textbf{45} (2003), 229--249.

\bibitem{Boffi+Rossi}
G.~Boffi and F.~Rossi. Lexicographic Gr\"obner bases for
transportation problems of format $r\times 3\times 3$. Preprint
submitted to Journal of Symbolic Computation.

\bibitem{Diaconis+Graham+Sturmfels}
P.~Diaconis, R.~Graham, and B.~Sturmfels. 
Primitive Partition Identities. {\sl Combinatorics--Paul Erd\"os is
eighty}, D.~Mikos, V.~Sos, T.~Szoni (eds.), 43-56. Bolyai Soociety
Mathematical Studies, 2, Budapest, 173--192.

\bibitem{Graver:75}
J.~E.~Graver. On the foundation of linear and integer programming
{I}. {\sl Mathematical Programming} \textbf{9} (1975), 207--226.

\bibitem{Hemmecke:PSP}
R.~Hemmecke.
On the positive sum property and the computation of {G}raver test sets.
{\sl Mathematical Programming} \textbf{96} (2003), 247--269.

\bibitem{Nairn:PhD-thesis}
K.~A.~Nairn. Graver Complexity of Monomial Curves in P3. PhD thesis,
Columbia University, 2003.


\bibitem{Santos+Sturmfels} F.~Santos and B.~Sturmfels. Higher
Lawrence configurations. J. Combin. Theory Ser. A {\bfseries 103}
(2003), 151--164.

\bibitem{Sturmfels+Thomas:97}
B.~Sturmfels and R.~R.~Thomas. Variation of cost functions in
integer programming. Mathematical Programming \textbf{77} (1997),
357--387.

\end{thebibliography}
\end{document}